\documentclass[12pt,a4paper]{article}
\usepackage[utf8]{inputenc}
\usepackage[english]{babel}
\usepackage[T1]{fontenc}
\usepackage{amsmath}
\usepackage{amsfonts}
\usepackage{bbm}
\usepackage{graphicx}
\usepackage{array}
\usepackage{authblk}
\usepackage[left=3cm,right=3cm,top=3cm,bottom=3cm]{geometry}

\usepackage{adjustbox}
\usepackage{transparent}

\usepackage{isomath}
\usepackage{mathtools}
\usepackage{amsthm}
\usepackage{slashed}
\usepackage{amssymb}
\usepackage{enumitem}
\usepackage{hyperref}
\hypersetup{colorlinks=true, 
citecolor=red}
\usepackage[new]{old-arrows}

\newtheorem{theorem}{Theorem}

\theoremstyle{definition}
\newtheorem{mdef}{{Definition}}[section]

\theoremstyle{definition}

\theoremstyle{definition}
\newtheorem{mrmk}{{Remark}}[section]

\theoremstyle{plain}
\newtheorem{mth}{Theorem}[section]

\theoremstyle{plain}
\newtheorem{mlem}{{Lemma}}[section]

\theoremstyle{plain}
\newtheorem{mprop}{{Proposition}}[section]

\theoremstyle{plain}
\newtheorem{mcor}{{Corollary}}[section]

\newcommand{\bigzero}{\mbox{\normalfont\Large 0}}
\newcommand{\rvline}{\hspace*{-\arraycolsep}\vline\hspace*{-\arraycolsep}}

\DeclarePairedDelimiter{\abs}{\lvert}{\rvert}
\makeatletter
\let\oldabs\abs
\def\abs{\@ifstar{\oldabs}{\oldabs*}}
\makeatother

\def\XXint#1#2#3{{\setbox0=\hbox{$#1{#2#3}{\int}$}
     \vcenter{\hbox{$#2#3$}}\kern-.5\wd0}}

\usepackage{tikz-cd} 
\usetikzlibrary{positioning, shapes, arrows.meta}
\tikzset{
    answer/.style={rectangle, draw, text width=15em, text badly centered, node distance=1cm, inner sep=0pt, minimum height=4em},
    block/.style={rectangle, draw, text width=10em, text centered},
     block2/.style={rectangle, draw, text width=5em, text centered},
     block3/.style={rectangle, draw, text width=5em, text centered, color=white},
}

\DeclareMathSymbol{\upLambda}{\mathalpha}{operators}{3}

\usepackage{pgfplots}
\usetikzlibrary{backgrounds}
\pgfplotsset{compat=1.14}
\usepackage{mathrsfs}
\usetikzlibrary{arrows}

\title{Statistical Fluctuation of Infinitesimal Spaces}
\author[1]{Damien Tageddine}
\author[1]{Jean-Christophe Nave}
\affil[1]{Department of Mathematics and Statistics, McGill University}
\date{}
\begin{document}
\maketitle
\begin{abstract}
This paper is a follow-up on the \emph{noncommutative differential geometry on infinitesimal spaces}  \cite{tageddine_noncommutative_2022}. In the present work, we extend the algebraic convergence from \cite{tageddine_noncommutative_2022}  to the geometric setting. On the one hand, we reformulate the definition of finite dimensional compatible Dirac operators using Clifford algebras. This definition also leads to a new construction of a Laplace operator. On the other hand, after a well-chosen Green function defined on a manifold, we show that when the Dirac operators can be interpreted as stochastic matrices. The sequence $(D_n)_{n\in \mathbb{N}}$ converges then in average to the usual Dirac operator on a spin manifold. The same conclusion can be drawn for the Laplace operator.
\end{abstract}
\newpage

{
  \hypersetup{linkcolor=black}
  \tableofcontents
}
\newpage
\section{Introduction}
The approximation theory of partial differential equations (PDE) can take various aspects. Traditionally, numerical analysis proposes different strategies to discretize operators. Depending on the situation, finite differences, finite elements, finite volumes, or spectral methods  may be utilized. In this process, the focus is usually analytical. That is, the aim is to control asymptotic convergence of the approximation error in a small parameter $(\Delta t, \Delta x,...)$. In fact, only a small subset of discretization techniques aim to preserve certain underlying structures (e.g. geometric, algebraic, etc…) of the continuous operator at the discrete level. The present work is a follow-up to \cite{tageddine_noncommutative_2022}, in which we provided a general framework  to tackle this question.\par 
In our paper, we are interested in the so-called family of compatible discretization, also called geometric discretization. The focus of these approximations is that discrete theory can, and indeed should, possess a geometric description on its own right. Among various types of approaches, one may mention the finite element exterior calculus \cite{arnold_compatible_2006}, the discrete exterior calculus \cite{desbrun_discrete_2005}, methods enforcing group symmetries such as in \cite{bihlo_invariant_2013,bihlo_convecting_2014}, or conservation laws such as those in \cite{wan_multiplier_2016, wan_arbitrarily_2018}. Since a fair amount of the theory of PDEs is developed on (subdomain of) $\mathbb{R}^n$, the latter approaches are, at least in their original incarnations, focused on Euclidean spaces. However, examples of partial differential equations arise in a wide variety of applications. As such, extensions of some of the previously mentioned techniques to non-Euclidean domains remains a challenge. Hence, a crucial question in the theory of discretization, is the generalization of classical geometric approaches to smooth manifolds. Also, and still on the topic of convergence analysis of finite elements, \cite{bobenko_convergence_2008} studies the cotangent discretization of the Laplace-Beltrami operator; the key result is that mean cruvature vectors converge in the sense of distributions, but fail to converge in $L^2$. Finally, there are the central research advances on diffusion maps in \cite{lafon_diffusion_2004,coifman_diffusion_2006} and the one on random point clouds in \cite{belkin_towards_nodate,belkin_constructing_2009}. As one can see, approaching the problem of compatible discretizations on manifolds rests heavily on the initial setup chosen to tackle it. The various results are therefore quite different, and perhaps appear disconnected from one another.\par 
In our recent work \cite{tageddine_noncommutative_2022}, we derive a general framework to describe finite difference calculus. This framework relies on the tools of $C^*$-algebras and noncommutative differential geometry (see \cite{connes_noncommutative_1994, connes_noncommutative_2007}). More specifically, starting uniquely from a differentiable manifold $M$, we exhibit a discrete space $X$ and its associate algebra $A(X)$ playing the role of an algebra of function. Then, using the natural setting of $C^*$-algebra, and its representation theory, we define a differential calculus on the space $X$. The corner stone of this definition is a so-called Dirac operator. Finally, this construction allows us to study spectral convergence with respect to a positive parameter $h$ in the case where $X$ is a lattice. Doing so may provide a general discrete construction of differential operators on smooth manifolds.\\
In the present work, we aim at studying the relation between the Dirac operator, $D$ defined in \cite{tageddine_noncommutative_2022} and its continuous counterpart $\mathcal{D}$, thus extending our previous construction beyond the case of a lattice.\par 
To this end, we recall the definition of the derivation $d$ given as a commutator with $D$:
\begin{equation}
da=\left[ D,a\right],
\end{equation}
where $a$ is an element of a $C^*$-algebra. This operator is analogous to the continuous Dirac operator $\mathcal{D}$. Indeed, in the case of a Riemannian manifold $(M,g)$ with a spinor bundle $\mathscr{S}\rightarrow M$, a \textit{Dirac operator} $\mathcal{D}$ on $\mathscr{S}$ is a differential operator whose principal symbol is that of $\mathsf{c}\circ d_{dR}$, where $\mathsf{c}$ is the quantization map. In particular, for any $a\in C^\infty(M)$, one has:
\begin{equation}
\left[ \mathcal{D},a\right] = i\mathsf{c}\circ d_{dR}(a).
\end{equation}\par
The operator $\mathcal{D}$ is defined using the Clifford algebras (see \cite[Def. 5.5.12 p.406]{rudolph_differential_2017} for a complete definition). Indeed, in local coordinates on a normal neighbourhood centred at a point $p$:
\begin{equation}
\mathcal{D}_p=\sum_{j=1}^de_j\left. \frac{\partial}{\partial x_j}\right|_{p}
\label{DiracOp}
\end{equation}
where $\lbrace e_j \ | \ j=1,\dots,d\rbrace$ is an orthonormal local frame embedded in a Clifford algebra. In \cite{tageddine_noncommutative_2022} , we gave the following definition for the operator $D$ :
\begin{equation}
D=\sum_{i<j}\omega_{ij}e_{ij}, \quad \omega_{ij}\in \mathbb{C}
\end{equation}
where the $e_{ij}$ are merely matrix elements in $M_{2n}(\mathbb{C})$ associated to a graph. In order to obtain an analogous description of \eqref{DiracOp} in terms of Clifford elements, in the present work, we redefine the finite dimensional Dirac operator as an element in the algebra $M_2(\mathbb{C})\otimes U(\mathfrak{g})$, where $U(\mathfrak{g})$ is the universal Lie algebra associated to a Lie algebra $\mathfrak{g}$. In addition, this definition leads to a construction of a Laplace operator using the inclusion map $U(\mathfrak{g})\rightarrow Cl(\mathfrak{g})$, where $Cl(\mathfrak{g})$ is the Clifford algebra on $\mathfrak{g}$.\par 
Now, to a fixed triangulation $X$, one can associate a collection of Dirac operators $(D_t)_{t\in \mathbb{N}}$, where each matrix $D_t$ can be seen as an irreducible matrix associate to the graph $G$ obtained from $X$. It has $n$ vertices labelled $1,\dots,n$, and there is an edge from vertex $i$ to a vertex $j$ precisely when $\omega_{ij}\neq 0$. More precisely, in the probabilistic setting, a vertex $i$ is connected to a vertex $j$ with probability $\omega_{ij}$.\\
Then, if we let $a_t\in \text{Dom}(D_t)$ and define the average operator, 
\begin{equation}
S_N=\frac{1}{N}\sum_{t=1}^Ne_t\left[  D_t,a_t\right]e^*_t,
\end{equation}
where $(e_t)_{t\in \mathbb{N}}$ is some family of projections. The key here is to choose the coefficients $\omega^t_{ij}$ associated to $D_t$ in order for the average operator $S_N$ to converge to $\left[ \mathcal{D} ,a\right]$ as $N\rightarrow \infty$. The main result of the present work is given by the following theorem.
\begin{mth}[Main result]
\label{main}
Let $\left\lbrace x^k_{i_0}\right\rbrace _{k=1}^n$ be a sequence of i.i.d. sampled points from a uniform distribution on an open normal neighbourhood $U_p$ of a point $p$ in a compact Riemannian manifold $M$ of dimension $d$. Let $\tilde{S}^{\hbar_n}_n$ be the associated operator given by:
\begin{equation}
\widehat{S}_n^{\hbar_n}:C^\infty(U_p)\rightarrow M_2(\mathbb{R})\otimes U(\mathfrak{gl}_{2m_n}), \qquad \widehat{S}_n^{\hbar_n} (a) :=\frac{1}{n}\sum_{k=1}^ne_k\left[ D^k_X, a_k\right]e_k^*.
\end{equation}
Put $\hbar_n=n^{-\alpha}$, where $\alpha>0$, then for $a\in C^\infty(U_p)$, in probability:
\begin{equation*}
\lim_{n\rightarrow\infty}\Psi\circ\widehat{S}_n^{\hbar_n}(a)=\left[ \mathcal{D},a\right](p).
\end{equation*} 
\end{mth}
Additionally, a similar result is proved for the case of the Laplace operator. These results generalize the previous ones obtained in \cite{tageddine_noncommutative_2022} beyond the case of a lattice.
It is worth mentioning at this point, in the realm of noncommutative geometry, the work of \cite{khalkhali_spectral_2021}, where in the same spirit (though in different context and approaches) the Dirac operators are defined as random matrices and form a so-called \textit{Dirac ensembles}. More specifically, the coefficients of the $N\times N$ matrix $D$ are random variables, following a prescribed density function; the spectral properties in the large $N$ limit are then explored. In the present work, however, the Dirac operator is associated to a graph and the density functions of the coefficients are specifically chosen to obtain a convergence result with respect to $\mathcal{D}$.\par
This paper is arranged as follows. In Section \ref{sect1}, we start with a presentation of the main results of \cite{tageddine_noncommutative_2022}. We then introduce Clifford algebras and universal enveloping algebras in order to define and study Dirac operators on finite dimensional spaces. In Section \ref{sect2}, we introduce a specific Hamilton-Jacobi equation with its associated Green function on $\mathbb{R}^d$ first, and their generalizations on a Riemannian manifold then. This is followed by some technical lemmas required to define the coefficients $\omega_{ij}$  necessary to prove Theorem \ref{main}. In Section \ref{sect3}, we prove our main result, Theorem \ref{main}, and we obtain as a by-product a convergence result for the Laplacian operator.
\section{Dirac operators in the algebraic setting}
\label{sect1}
In this section, we introduce two of the main algebraic tools that we are going to use in this study: the Clifford algebras and the universal enveloping Lie algebra. We then define and study Dirac operators on finite spectral triples in terms of root vectors of a Lie algebra $\mathfrak{g}$.
\subsection{Noncommutative Geometry on Infinitesimal spaces}
In the research paper \cite{tageddine_noncommutative_2022}, we show that a discrete topological space $X$ can be identified to the spectrum $Spec(A)$ of a $C^*$-algebra $A$. Starting with a Riemannian manifold $(M,g)$, we construct an inverse system of triangulations, $(K_n)$ which become sufficiently fine for large $n$. Using the Behncke-Leptin construction, we associate to each $K_n$ a $C^*$-algebra $A_n$ such that the triangulation $K_n$ is identified with the spectrum $Spec(A_n)$. We then form an inductive system $(A_n)$ with limit $A_\infty$.
\begin{theorem}
The spectrum $Spec(A_\infty)$ equipped with the hull-kernel topology is homeomorphic to the space $X_\infty$ and 
\begin{equation}
\lim_{\leftarrow}Spec(A_i) \simeq Spec(\lim_{\rightarrow}A_i).
\end{equation}
\end{theorem}
We then show that the centre of $A_\infty$ is isomorphic to the space of continuous function $C(M)$. In this sense, any element $g\in C(M)$ can be uniformly approximated arbitrarily closely by elements $a_n$ in the central subalgebras $\mathfrak{A}_n$.
\begin{theorem} The space of continuous function $C(M)$ is approximated by the system of commutative subalgebras $(\mathfrak{A}_n,\phi^*_{n,\infty})$ in the following sense:
\begin{equation}
C(M)=\overline{\bigcup_{n\in \mathbb{N}}\phi^*_{n,\infty}(\mathfrak{A}_n)} \cap C(M).
\end{equation}
\end{theorem}
In addition, the sequence of representations $(H_n)$ is also considered as a direct system with limit $H_\infty$ containing the space of square integrable functions $L^2(M)$. Finally, we define the spectral triples $(\mathfrak{A},\mathfrak{h},D_n)$, where $D_n$ is a so-called Dirac operator. We show that under certain conditions, the sequence $(D_n)$ converges to the multiplication operator by the de Rham differential $d_ca$.
\begin{theorem} (Spectral convergence)
There exists a finite measure $\mu$ and a unitary operator 
\begin{equation}
U:L^2(\mathbb{R})\rightarrow L^2(\mathbb{R},d\mu)
\end{equation}
such that, 
\begin{equation}
U[D,a]U^{-1}\phi=\frac{d a}{d x}\phi, \quad \forall \phi \in L^2(\mathbb{R}),
\end{equation}
Moreover, the norm of the commutator is given by $\|\left[ D,a\right] \|=\|d_ca\|_\infty$.
\label{spectral}
\end{theorem}
Thus, we have built a correspondence between a given triangulation $X$ and a Dirac operator $D$: the non-zero coefficients of $D$ are determined by the connectivity between vertices of the graph. We showed that this is however not enough to represent the metric of the manifold. Thus, we ask now the question on how to set the coefficients $\omega_{ij}$ of $D$ so that at the limit (in the sense of \eqref{commutator}) the sequence converges.
\subsection{Clifford algebras}
Let $V$ be a finite dimensional vector space over a commutative field $\mathbb{K}$ of characteristic zero endowed with a quadratic form $q$. Let $T(V)$ be the tensor algebra over $V$. Consider the ideal $I_q$ in $T(V)$ generated by all elements of the form $v\otimes v + q(v)$ for $v\in V$. Then the quotient algebra 
\begin{equation}
Cl(V,q)=T(V)/I_q.
\end{equation}
is the Clifford algebra associated to the quadratic space $(V,q)$.\\
Moreover, we can choose any orthonormal basis $Z_i$ of $V$ with respect to $q$ as a set of generators of $Cl(V)$. We then have the relations, 
\begin{equation}
Z_iZ_j=-Z_jZ_i, \quad i\neq j, \quad Z_i^2=-1.
\end{equation}
Then the following set 
\begin{equation}
Z_{i_1}Z_{i_2}\cdots Z_{i_k} \quad 1\leq <i_1<i_2<\cdots < i_k\leq n=\dim V
\end{equation}
spans $Cl(V)$. In addition, given a $q$-orthonormal basis $Z_i$ of $V$, the mapping 
\begin{equation}
1\mapsto 1, \quad Z_{i_1}\cdots Z_{i_k}\mapsto Z_{i_1}\wedge \cdots \wedge Z_{i_k}
\end{equation}
yields an isomorphism of vector spaces $Cl(V,q)\simeq \bigwedge V$.
\subsection{Dirac operators in the Clifford algebra setting}
Let $k=\mathbb{R}, \mathbb{C}$ and let $\mathfrak{g}$ be a Lie algebra over $k$. We start by recalling the definition of the universal enveloping algebra.
\begin{mdef}
The $\textit{universal enveloping algebra}$ of $\mathfrak{g}$ is a map $\varphi:\mathfrak{g}\rightarrow U(\mathfrak{g})$, where $U(\mathfrak{g})$ us a unital associative algebra, satisfying the following properties:
\begin{itemize}
\item[1)] $\varphi$ is a Lie algebra homomorphism, i.e. $\varphi$ is $k$-linear and
\begin{equation}
\varphi\left( \left[X,Y \right] \right) =\varphi(X)\varphi(Y) - \varphi(Y)\varphi(X).
\end{equation}
\item[2)] If $A$ is any associative algebra with a unit and $\alpha:\mathfrak{g}\rightarrow A$ is any Lie algebra homomorphism, there is a unique homomorphism of associative algebras $\beta:U(\mathfrak{g})\rightarrow A$ such that the diagram 
\begin{center}
\begin{tikzcd}
  \mathfrak{g} \arrow[r,"\varphi"] \arrow[d,"\alpha"']
    & U(\mathfrak{g}) \arrow[dl,"\beta"]\\
	A 
\end{tikzcd}
\end{center}
is commutative, i.e. there is an isomorphism 
\begin{equation}
\text{Hom}_{Lie}(\mathfrak{g}, LA)\simeq \text{Hom}_{Ass}(U(\mathfrak{g}),A)
\end{equation}
\end{itemize}
\end{mdef}
We will now give a definition of Dirac operators on finite spectral triples in terms of root vectors of a Lie algebra $\mathfrak{g}$. Then, using the canonical embedding $\mathfrak{g}\hookrightarrow Cl(\mathfrak{g})$ into the Clifford algebra, we define a Laplace-type operator.\par
Consider the algebra $A=\mathfrak{gl}_{2N}(\mathbb{C})$ of complex matrices with its standard Lie algebra structure. In \cite{tageddine_noncommutative_2022}, we have introduced the finite dimensional spectral triple $(\mathfrak{A},\mathfrak{h},D)$ given by:
\begin{itemize}
\item[•]$\mathfrak{A} $ is a Cartan subalgebra of Lie subalgebra $\mathfrak{g}$ of $A$,
\item[•]$\mathfrak{h}=\mathbb{C}^{2N}$,
\item[•]$\gamma = \left( \begin{array}{cc}
1_N & 0 \\ 
0 & -1_N
\end{array} \right) $.
\end{itemize}
The chirality element $\gamma$ induces a decomposition of the representation space $\mathfrak{h}$ into the eigenspaces $\mathfrak{h}^{\pm}$ corresponding to the eigenvalues $1$ and $-1$ such that $\mathfrak{h}=\mathfrak{h}^+\oplus \mathfrak{h}^{-}$. Incidentally, one has the decomposition of the algebra $A$ as follows:
\begin{equation}
\mathfrak{gl}_{2N}=\mathfrak{gl}_{2N}^+\oplus \mathfrak{gl}_{2N}^{-}.
\label{decomp}
\end{equation}
Notice then that the pair $(\mathfrak{gl}_{2N}^{+}, \mathfrak{gl}_{2N}^{-})$ forms a Cartan pair.
Any endomorphism $a\in \text{End}(\mathfrak{h})$ defines an endomorphism $\rho_a\in \mathfrak{gl}_{2N}^+$ given by
\begin{equation}
\rho_a=\left( \begin{array}{cc}
a & 0 \\ 
0 & -a^T
\end{array} 
\right)\in \mathfrak{sp}(2N,\mathbb{C})\cap \mathfrak{gl}_{2N}^+.
\end{equation}
We consider the compact real case with the embedding 
\begin{equation}
\mathfrak{sp}(N)= \mathfrak{sp}(2N,\mathbb{C})\cap \mathfrak{u}(2N)\hookrightarrow \mathfrak{so}(4N).
\end{equation}
If $a$ is a diagonal element of $\text{End}(\mathfrak{h}^+)$, the map $a\mapsto \rho_a$ identifies $a$ with an element of the maximal commutative subalgebra $\mathfrak{t}$ of $\mathfrak{so}(4N)$:
\begin{equation}
\mathfrak{t}=\left\lbrace 
\left( 
\begin{array}{ccc}
A_1 & 0 & 0 \\ 
0 & \ddots & 0 \\ 
0 & 0 & A_n
\end{array} 
\right) , A_j=
\left( 
\begin{array}{cc}
0 & a_j \\ 
-a_j & 0
\end{array} 
\right) 
\right\rbrace 
\end{equation}
Consider the Cartan subalgebra $\mathfrak{A}=\mathfrak{t}+i\mathfrak{t}$ of $\mathfrak{so}(4N,\mathbb{C})$. The root vectors are $4N\times 4N$ block matrices having $2\times 2$-matrix $C_s$, $s\in \left\lbrace 1,\dots, 4\right\rbrace $
\begin{equation}
X=\left( 
\begin{array}{cc}
0 & C_s \\ 
-C^t_s & 0
\end{array} 
\right) 
\end{equation}
in the position $(i,j)$ with $i<j$ and where 
\begin{align*}
&C_1=\left( 
\begin{array}{cc}
1 & i \\ 
i & -1
\end{array} 
\right), \quad
C_2=
\left( 
\begin{array}{cc}
1 & -i \\ 
-i & -1
\end{array} 
\right) ,\quad C_3=\left( 
\begin{array}{cc}
1 & -i \\ 
i & 1
\end{array} 
\right), \quad
C_4=
\left( 
\begin{array}{cc}
1 & -i \\ 
i & -1
\end{array} 
\right).
\end{align*}
associated to the linear functional in $\mathfrak{H}^*$ given by $i(a_i+a_j)$, $-i(a_i+a_j)$, $i(a_i-a_j)$ and $i(a_j-a_i)$.\par 
We will denote by $\mathfrak{g}$ the Lie algebra $\mathfrak{so}_{4N}$. We then consider the unital associative algebra $M_2(\mathbb{C})\otimes \mathfrak{gl}_{2N}$ and the homomorphism:
\begin{equation}
\varphi:\mathfrak{g}\rightarrow M_2(\mathbb{C})\otimes \mathfrak{gl}_{2N},\quad
\varphi(X)=\sum_{1\leq i,j\leq 2N}X_{ij}\otimes E_{ij},
\label{hom1}
\end{equation}
where $E_{ij}$ is the standard basis in $\mathfrak{gl}_{2N}$ and $X_{ij}$ are the $2\times 2$-submatrix of $X=(x_{rs})$ obtained by keeping $i+1\leq r\leq i+2$ and $j+1\leq s\leq j+2$. In addition, $\varphi $ is a Lie algebra homomorphism with
$\varphi(\left[X, Y\right] )=\varphi(X)\varphi(Y)-\varphi(Y)\varphi(X)$.\\
Then, using the universal property of $U(\mathfrak{g})$, the map $\varphi$ extends into the homomorphism $\widehat{\varphi}:U(\mathfrak{g})\rightarrow M_2(\mathbb{C})\otimes \mathfrak{gl}_{2N}$. Furthermore, taking the canonical embedding $h:\mathfrak{gl}_{2N}\rightarrow U(\mathfrak{gl}_{2N})$, we get by composing the Lie algebra homomorphism 
\begin{equation}
h\circ \widehat{\varphi}:U(\mathfrak{g})\rightarrow M_2(\mathbb{C})\otimes U(\mathfrak{gl}_{2N}).
\label{hom2}
\end{equation} 
Let $\left\lbrace Z_{ij} \right\rbrace $ be an orthonormal basis of root vectors in $\mathfrak{g}$, associated to the root $-i(a_j+a_k)$, we define the operator $W$ by
\begin{equation}
W=\sum_{i,j}\omega^W_{ij}Z_{ij}
\label{opA}
\end{equation}
as an element of $U(\mathfrak{g})$, where $\omega^W_{ij}$ are real coefficients.
\begin{mdef} Given an operator $W$ as in \eqref{opA}, a Dirac operator $D_W$ is an element of $M_2(\mathbb{C})\otimes U(\mathfrak{gl}_{2N})$ defined by:
\begin{equation}
D_W=\frac{i}{\hbar}\mathfrak{Re}(W),
\label{DiracW}
\end{equation}
where $\hbar>0$ is a real parameter.
\end{mdef}
\begin{mrmk}
In the previous definition, $D_W$ depends on the choice of element $W$ and in fact, more specifically on the choices of basis elements $Z_{ij}$. Another definition, independent on the choice of basis elements, of Dirac operators on Lie algebras can be found in \cite{meinrenken_clifford_2013} 
\end{mrmk}
\begin{mlem} Let $C_2=X+iY$ be the root vector associated to the root $-i(a_i+a_j)$. Fix an element $W$ as in \eqref{opA}. Then, for any $a\in \mathfrak{A}$, the exterior derivative can be written as:
\begin{equation}
\left[ D_W,a\right] = \frac{i}{\hbar}\sum_{i,j}\omega^W_{ij}\alpha_{ij}(a)Y\otimes E_{ij},
\label{dirac}
\end{equation}
an element of $M_2(\mathbb{C})\otimes U(\mathfrak{gl}_{2N})$ and with $\alpha_{ij}=a_i-a_j$.
\label{LemDA}
\end{mlem}
\begin{proof}
From the definition of $D_W$ and the definition of root vectors, we get that:
\begin{equation}
\left[ D_W, a\right] =\frac{i}{2\hbar}\sum_{i,j}\omega^W_{ij}(a_i-a_j)Z_{ij} - \frac{i}{2\hbar}\sum_{i,j}\omega^W_{ij}(a_i-a_j)Z_{ij}^*.
\end{equation}
Then, using the map $h\circ \widehat{\varphi}$, given by \eqref{hom1} and \eqref{hom2}, we can identify a basis element $Z_{ij}$ with an element in $M_2(\mathbb{C})\otimes U(\mathfrak{g})$ of the form $C_2\otimes E_{ij}$. Hence, we have that:
\begin{equation}
\left[ D_W, a\right] =\frac{i}{2\hbar}\sum_{i,j}\omega^W_{ij}(a_i-a_j)C_2\otimes E_{ij} - \frac{i}{2\hbar}\sum_{i,j}\omega^W_{ij}(a_i-a_j)C_2^*\otimes E_{ij}^t.
\end{equation}
Simplifying this expression using the fact that $E_{ij}^t=E_{ji}$, we get:
\begin{equation}
\left[ D_W, a\right] =\frac{i}{\hbar}\sum_{ij}\omega^W_{ij}\alpha_{ij}(a)Y\otimes E_{ij}.
\end{equation}
with $\alpha_{ij}(a)=a_i-a_j$.
\end{proof}
Furthermore, we recall that there exists a canonical Lie algebra homomorphism $\psi:\mathfrak{gl}_{2N}\rightarrow Cl(\mathfrak{gl}_{2N})$ which extends into the map on the universal enveloping Lie algebra:
\begin{equation}
\widehat{\psi}:U(\mathfrak{gl}_{2N})\rightarrow Cl(\mathfrak{gl}_{2N}). 
\label{mapx}
\end{equation}
We use this map to define a Laplace operator.
\begin{mdef}[Laplacian] Fix an element $W$. We then define the Laplace operator $\Delta$ on $\mathfrak{A}$ using the non-graded commutator. For any $a\in \mathfrak{A}$
\begin{equation}
\Delta(a):=\frac{1}{2}\widehat{\psi}(\left[ D_W,\left[D_W,a \right] \right] ).
\label{Laplacianop}
\end{equation}
\end{mdef}
\begin{mprop}
\label{laplace}
For any $a\in \mathfrak{A}$, the Laplace operator is given by
\begin{equation}
\Delta(a)=-\Omega_\mathfrak{g}(a)\otimes 1
\end{equation}
where $\Omega_\mathfrak{g}=\frac{1}{\hbar^2}\sum_{i,j}\omega_{ij}^2J\otimes\alpha_{ij}$ is an element of $End(\mathfrak{A},M_2(\mathbb{C}))$.
\end{mprop}
\begin{proof}
Let $D_W$ be a Dirac operator, then the bi-commutator of the Laplacian gives::
\begin{equation}
\left[ D_W, \left[ D_W, a\right]\right] =D_W \left[ D_W, a\right] - \left[ D_W, a\right]D_W.
\end{equation}
Thus, using Lemma \ref{LemDA}, we obtain
\begin{equation}
\left[ D_W, \left[ D_W, a\right]\right]=\frac{2}{\hbar^2}\sum_{ij}\omega_{ij}^2\alpha_{ij}(a)J\otimes E_{ij}^2  +\frac{1}{\hbar^2}\sum_{(ij)\neq (kl)}\omega_{ij}\omega_{kl}\alpha_{kl}(a)J \otimes \left[ E_{ij}, E_{kl}\right]_+
\label{bicom}
\end{equation}
with the bracket $[A,B]_{+}=AB+BA$ and where the matrix $J$ is given by:
\begin{equation*}
J=\left( 
\begin{array}{cc}
0 & -1 \\ 
1 & 0
\end{array} 
\right).
\end{equation*}
Finally, applying the map $\widehat{\psi}$, the second term of the left-hand-side in Equation \eqref{bicom} vanishes and we get: 
\begin{equation*}
\Delta(a)=-\frac{1}{\hbar^2}\sum_{i,j}\omega_{ij}^2\alpha_{ij}(a)J\otimes 1.
\end{equation*}
\end{proof}
We have kept the definition of the Dirac operator $D_W$ in \eqref{DiracW} very general, however we recall that the operator we are interested in are the \textit{compatible} ones with respect to a graph $X$. In other words, the value $\omega_{ij}$ is non-zero if $ij$ is an edge in $X$.\par
Now, let us recall how the space $X$ is obtained from a manifold $M$; more details can be found in \cite{tageddine_noncommutative_2022}. One starts with a triangulation of $M$ and then consider the \textit{dual of the triangulation} that we will call $X$. In fact, in \cite{tageddine_noncommutative_2022} we used a slightly different terminology and considered the triangulation as a poset, then looked at the opposite poset with reversed order.\par
Since we are working with a graph $X$ obtained from a dual triangulation, every vertex $i$ has exactly $d+1$ neighbours, i.e. only $d+1$ of the $\omega_{ij}$ are non-zero for a fixed $i$. Hence, if we fix a vertex $i_0$, the definition of the commutator with $D_W$ becomes:
\begin{equation}
\left( \left[ D_W,a\right]\right)_{i_0} = \frac{i}{\hbar}\sum_{j=1}^{d+1}\omega^W_{i_0k_j}\alpha_{i_0k_j}(a)Y\otimes E_{i_0k_j},
\end{equation}
Here, we relabel the index $j$ without lost of generalities and to keep this indexing simple. We will also drop the $W$ index for the same reasons and get:
\begin{equation}
\left( \left[ D,a\right]\right)_{i_0} = \frac{i}{\hbar}\sum_{j=1}^{d+1}\omega_{i_0j}\alpha_{i_0j}(a)Y\otimes E_{i_0j},
\label{dirac2}
\end{equation}
Finally, let us recall the (true) Dirac operator on a manifold is given in local coordinates on a normal neighbourhood centred at a point $p$:
\begin{equation}
\mathcal{D}_p=\sum_{j=1}^de_j\left. \frac{\partial}{\partial x_j}\right|_{p}
\label{DiracOp}
\end{equation}
where $\lbrace e_j \ | \ j=1,\dots,d\rbrace$ is an orthonormal local frame embedded in the Clifford algebra $Cl(\mathbb{R}^d)$ using the natural embedding $\mathbb{R}^d \subset Cl(\mathbb{R}^d)$.\par
Nevertheless, the Dirac operator as expressed in \eqref{dirac2} is not an element of a Clifford algebra. Moreover, the dimensions do not match. Indeed, because of the structure of the triangulation, there are $d+1$ independent vectors in the expression \eqref{dirac2}, instead of $d$ as the dimension of the manifold $M$. Since we are trying to approximate the true Dirac operator in \eqref{DiracOp}, we need to re-write Equation \eqref{dirac2} in terms of Clifford elements in dimension $d$. To do so, let us denote by $V_{i_0}$, the vector space defined by:
\begin{equation}
V_{i_0}:=\text{span}\left\lbrace Y\otimes E_{i_0,1}, \cdots , Y\otimes E_{i_0,d+1} \right\rbrace .
\end{equation}
Then, consider the isomorphism:
\begin{equation}
\tau:V_{i_0}\xrightarrow[]{\simeq}\mathbb{R}^{d+1} \quad \tau (Y\otimes E_{i_0,j})=\widehat{e}_j, \ \forall 1\leq j\leq d+1
\end{equation}
where $\left\lbrace \widehat{e}_j\right\rbrace_{j=1}^{d+1} $ is the canonical basis on $\mathbb{R}^{d+1}$ with respect to the standard inner product. Moreover, defines the projection $p$ on the subspace spanned by $\left\lbrace \widehat{e}_j\right\rbrace_{j=1}^{d}$ and identified with $\mathbb{R}^d$. Finally, if let the embedding $\rho : \mathbb{R}^d\rightarrow Cl(\mathbb{R}^d)$, we can compose these maps and define:
\begin{equation}
\Psi:=\rho \circ p\circ \tau: V_{i_0}\rightarrow Cl(\mathbb{R}^d), \quad \Psi\left( \left[ D,a\right]_{i_0}\right) =  \frac{i}{\hbar}\sum_{j=1}^{d}\omega_{i_0j}\alpha_{i_0j}(a)e_j
\label{mapPsi}
\end{equation}
which allows us to express the commutator in terms of Clifford elements $e_j$. We notice, nevertheless, that this construction is not canonical and depends on the choice of isomorphism $\tau$. 
\subsection{Perron-Frobenius bound on $\left[D,a \right] $}
To conclude this section, and before being able to show a convergence result to the Dirac operator $\mathcal{D}$, we would like to prove a preliminary result on the commutator $[D,a]$ and its boundedness at the limit when $\hbar\rightarrow 0$. This result follows from the correspondence between $D$ and the graph associated, using the Perron-Frobenius theorem. We only need to consider the operator $D$ as a compatible operator in some matrix space, without relying on the Clifford algebra setting\par
We consider an infinite collection $\left\lbrace \mathfrak{A}_n : n\in \mathbb{N} \right\rbrace $ of commutative $C^*$-algebras. In this case, we have identified each of the $\mathfrak{A}_n$ with the Cartan subalgebras $\mathrm{h}_i$ inside the finite dimensional algebras $B_n=\mathfrak{so}_{2m_n}(\mathbb{C})$ where $m_n\rightarrow \infty$ when $n\rightarrow \infty$. We can then construct the product: 
\begin{equation}
B_\omega= \prod_{n\in \mathbb{N}}B_n=\lbrace (a_n) : \|a_n\|=\sup \|a_n\|<\infty\rbrace .
\end{equation}
Let $a$ be an element in $C^\infty(M)$, then there exists a coherent sequence 
$(a_i)$ such that
\begin{equation}
a=(a_0, a_1, \cdots, a_n, \cdots)\in  \prod_{n\in \mathbb{N}} \mathfrak{A}_n.
\end{equation}
We define a spectral triple on $B_\omega$ by introducing the limit Dirac operator $D$ as the sequence 
\begin{equation}
D=(D_0,D_1,\cdots, D_n, \cdots ) \in \prod_{n\in \mathbb{N}} \mathfrak{gl}_{2m_n}^-(\mathbb{C}),
\end{equation}
where each $D_i$ is a Dirac operator associated to a poset $X_i^{op}$ in the sense of \cite{tageddine_noncommutative_2022}.
This in turns induces a spectral triple on $ \prod_{n\in \mathbb{N}} \mathfrak{A}_n$ along with the commutator:
\begin{equation}
d_Da:=[D,a]=([D_0,a_0],[D_1,a_1],\cdots, [D_n,a_n], \cdots )\in \prod_{n\in \mathbb{N}} \mathfrak{gl}_{2m_n}^-(\mathbb{C}).
\label{commutator}
\end{equation}
In order to show that $[D,a]$ is a bounded operator, we use Perron-Frobenius theorem, which we start by recalling.
\begin{mth}[Perron-Frobenius]
Let $A=(a_{ij})$ be an $n\times n$ positive matrix: $a_{ij}>0$ for $1\leq i, j\leq n$. Then there exists a positive real number $r$, called the \textit{Perron-Frobenius eigenvalue}, such that $r$ is an eigenvalue of $A$. Moreover, if the spectral radius $\rho(A)$ is equal to $r$. \\
The Perron-Frobenius eigenvalue satisfies the inequalities:
\begin{equation*}
\min_{i}\sum_{j}a_{ij}\leq  r \leq \max_{i}\sum_{j}a_{ij}.
\end{equation*}
\label{Per-Frob}
\end{mth}
\begin{mprop} For any $a\in \mathfrak{A}$, the spectral radius $\rho(d_Da)$ of $d_Da$ is bounded by  
\begin{equation}
\rho(d_Da)\leq \|d_{dR}a\|_{\infty}.
\end{equation}
\end{mprop}
\begin{proof}
We consider the sequence of Dirac operators $(D_\alpha)_{\alpha\in \mathbb{N}}$ associated to $D$. \par
Let $\varepsilon>0$ and $\alpha\in \mathbb{N}$ and define the operator $\widetilde{d_{D_\alpha}a}$ such that
\begin{equation}
(\widetilde{d_Da})_{ij}=\left\lbrace  \begin{array}{cc}
|(d_{D_\alpha}a)_{ij}| & \text{if $(d_{D_\alpha}a)_{ij} \neq 0$} \\ 
\varepsilon  & \text{otherwise}
\end{array} 
\right. 
\end{equation}
The matrix $\widetilde{d_{D_\alpha}a}$ is positive by construction. In addition, we have the upper-bound:
\begin{equation}
\|\left( d_{D_\alpha}a\right) ^k\|^2_F\leq \|(\widetilde{d_{D_\alpha}a})^k\|^2_F.
\end{equation}
Hence, using Theorem \ref{Per-Frob}, we deduce that
\begin{align*}
\rho(d_{D_\alpha}a)^2=\lim_{k\infty}\|(d_{D_\alpha}a)^k\|_F^{\frac{2}{k}}&\leq \lim_{k\infty}\|(\widetilde{d_{D_\alpha}a})^k\|_F^{\frac{2}{k}}\\
&=\rho(\widetilde{d_{D_\alpha}a})^2\\
&\lesssim\max_{1\leq i\leq n} \sum_{j} |(d_{D_\alpha}a)_{ij}|^2 +N\varepsilon^2.
\end{align*}
The value of $N$ is the number of nonzero coefficient in $(d_{D_\alpha}a)_{ij}$ and thus only depends on the number of adjacency vertex in $X^{op}_\alpha$ which by definition equal to $d+1$, where $d$ is the dimension.\\
Hence there exists a positive constant $C_M$, which depends on the maximal length of geodesics ($M$ is compact) but is independent of $\alpha$, such that 
\begin{equation}
\rho(d_{D_\alpha}a)^2\leq C_M \|d_{dR}a\|_{\infty}^2 + (d+1)\varepsilon^2.
\end{equation}
The last inequality holds for an arbitrary $\varepsilon >0$ and $\alpha\in \mathbb{N}$. The result follows  then by taking $\varepsilon$ to $0$. 
\end{proof}
\begin{mcor} For each $a\in \mathfrak{A}$, the operator $\left[ D,a\right] $ is a bounded operator.
\end{mcor}
\begin{mrmk} It is clear that in the following framework, not only the Dirac operator $D$ define a differential structure, but it also plays the role of a transition matrix. This last point will be made clearer in the following section.
\end{mrmk}
\section{Green's function and integral operators}
\label{sect2}
In this section, we are going to introduce the Green's function of a suitable second-order Cauchy problem. From it, we exhibit a probability distribution that will be used in Section \ref{sect3} in the definition of the Dirac operator. Furthermore, we will extend this probability distribution defined in $\mathbb{R}^d$ to a manifold using the exponential map and prove some technical lemmas.\par
\subsection{Hamilton-Jacobi equation with vanishing viscosity}
For reasons that will be detailed in this section, we are interested in the following Hamilton-Jacobi equation with vanishing viscosity over the Euclidean space $\mathbb{R}^d$:
\begin{equation}
\begin{aligned}
&\partial_t u = s\cdot \nabla u + \varphi(t) \Delta u\\
& u(x,0)=u_0(x) 
\end{aligned}
\label{HJeqt}
\end{equation}
where $s$ is a unit vector in $\mathbb{R}^d$, $u_0$ is a smooth initial condition and $\varphi\in C(\mathbb{R}_+)$ satisfying:
 \begin{equation}
\varphi(t) = O_{t\rightarrow 0^+}(t).
 \end{equation}
 The fundamental solution of Equation \eqref{HJeqt} is obtained by taking the initial condition to be the Dirac distribution $\delta(y-x)$ for $y\in \mathbb{R}^d$ fixed. The normalized fundamental solution denoted by $G$ is given by:
\begin{equation}
G_t(x,y)=\frac{1}{(4\pi \Phi(t))^{\frac{d}{2}}}\exp\left(-\frac{|y-x+st|^2}{4\Phi(t)}\right)dx, \qquad \text{with \ $\Phi(t)=\int_0^t\varphi (s)ds$}.
\label{Green}
\end{equation}
Hence, the general solution can be obtained by convolution:
\begin{equation}
u(y,t)=\left(u_0*G\right)(y,t)
\end{equation}
then we have the following lemma.
\begin{mlem}  Let $u$ be the solution of Equation\eqref{HJeqt} with initial condition $u_0\in \mathbb{C}^\infty(\mathbb{R}^d)$, then it satisfies the initial condition:
\begin{equation}
\partial_tu|_{t=0}=s\cdot \nabla u_0.
\end{equation}
\label{FP-th}
\end{mlem}
\begin{proof}
First, let us notice that the family $\left\lbrace G_t \right\rbrace_{t>0}$ is an approximation of the identity:
\begin{equation}
\forall f\in \mathbb{C}^\infty(\mathbb{R}^d), \quad \lim_{t\rightarrow 0}\|f-f*G_t\|_\infty.
\end{equation}
In addition, since the differential operator in Equation \eqref{HJeqt} has constant coefficients in the $y$ variable, we have:
\begin{equation}
\partial_tu=(s\cdot \nabla  + \varphi \Delta )u_0*G
\end{equation}
from which we get by taking the limit when $t$ goes to $0$:
\begin{equation}
\partial_tu(y,0)=s\cdot \nabla u_0(y).
\end{equation}
\end{proof}
Therefore, we can define the probability measure $d\nu$ given by 
\begin{equation}
d\mu_{y,t}(x)=G_t(x,y)dx 
\label{density}
\end{equation}
with respect to the Lebesgue measure on $\mathbb{R}^d$. Hence, we define the distribution:
\begin{equation}
T_{\mu_{y,t}}(\varphi)(y,t)=\int_{\mathbb{R}^d}\varphi(x) d\mu_{y,t}(x),
\end{equation}
that will focus our attention in the next section.
\begin{mlem}[Reduction to a ball] Consider an open ball $B_\delta\subset \mathbb{R}^n$ of radius $delta>0$ such that $p\in B_\delta$. For any function $f\in L^\infty(B_\delta)$ and a smooth extension $\overline{f}$ of $f$ in $L^\infty(\mathbb{R}^d)$ , we have, as $t\rightarrow 0$:
\begin{equation}
\left|\int_{B_\delta} G_t(x,y)f(x)dx-\int_{\mathbb{R}^d} G_t(x,y)\overline{f}(x)dx\right| = o(t^{d})
\end{equation}
\label{Reduc}
\end{mlem}
\begin{proof}
Without lost of generality, we may take $y=0$; after a change of variable $u=\frac{x}{2\sqrt{\Phi(t)}}$, we see that:
\begin{equation*}
\left|\int_{B_\delta} G_t(x,y)f(x)dx-\int_{\mathbb{R}^d} G_t(x,y)f(x)dx\right| \leq \frac{\|\overline{f}\|_\infty\Phi(t)^{\frac{d}{2}}}{\pi^{\frac{d}{2}}}e^{-\frac{t^2}{4\sqrt{\Phi(t)}}}\int_{B_{\delta_t}^c}e^{-\frac{\delta}{4}|u|}dx
\end{equation*}
with $\delta_t=\frac{\delta}{2\Phi(t)}$ and where $\Phi(t)$ decreases as $t$ tends to zero.
\end{proof}
Thus, we may equivalently consider the distribution $T_\mu$ restricted on a ball $B_\delta$.
\subsection{Some remarks families of one-parameter operators}
The previous construction of the Green function of the Cauchy problem \eqref{HJeqt} can be reformulated in the more general setting of families of one-parameter operators and the so-called Fokker-Planck equation \cite{bogachev_fokkerplanckkolmogorov_2015}. The construction goes as follows.\par
 Let $L^2(M)$ considered as a Banach space for the Lebesgue measure $\nu$.  Let then $U$ be an open subset of $M \times \mathbb{R}_+ $. We denote by $\mathscr{D}(U)$ the set of \textit{test functions} on $U$ and its topological dual $\mathscr{D}'(U)$ the space of \textit{distribution}. A family of one-parameter operators $\lbrace P_t\rbrace_{t\geq 0}$ is a family of linear operators on $L^2(M)$ defined by:
\begin{equation}
P_0=id, \qquad (P_tf)(x)=\int_{M}f(y)p_t(x,y)d\nu(y)
\end{equation}
such that $p_t(x,y)$ is a $\nu \times \nu$-measurable function on $M\times M$. Now let us define the new measure $\mu_{t,x}$, for $t\geq 0$ and $x\in M$, by 
\begin{equation}
\mu_{x,t}(A)=\int_{A}p_t(x,y)d\nu(y),
\end{equation}
for any $\nu$-measurable subset $A$. Assume that $\mu_{t,x}$ is a probability measure for every $(x,t)\in M\times \mathbb{R}_+ $. In addition, we assume that $P_t$ admits a derivative $\partial_t$. Then, one can associate to any operator $P_t$ a distribution; let $\varphi\in \mathscr{D}(U)$ and define 
\begin{equation}
    T_{\mu_{x,t}}(f)=\int_M f(y)d\mu_{x,t}(y)
\end{equation}
as a map on $\mathscr{D}'(U)$. 
In the special case where
\begin{equation}
\lim_{t\rightarrow 0+}\frac{\partial^k\widehat{\mu}_{x,t}(0)}{t}=0, \quad \forall k\geq 3,
\label{limite}
\end{equation}
then $\mu_{x,t}$ satisfy the parabolic equation
\begin{equation}
\left. \frac{\partial \mu_{x,t}}{\partial t}\right|_{t=0} = \mathscr{L}_{A,b}(\mu_{x,0})
\label{FP-eq}
\end{equation}
in the weak sense, called the \textit{Fokker-Planck equation}. The operator $\mathscr{L}_{A,b}$ is given by
\begin{equation}
\mathscr{L}_{A,b}f=tr(AD^2f)+\left\langle b,\nabla f \right\rangle, \quad f\in C^\infty_c(M)
\end{equation}
and where $A=(a^{ij})$ is a mapping on $M$ with values in the space of nonnegative symmetric linear operator on $\mathbb{R}^d$ and $b=(b^i)$ is a vector field on $M$. In the special case of Equation \eqref{HJeqt}, we have $A=1$ and $b=s$.\par
Following this idea, it is interesting to consider semigroup machinery as another approach to the problem of approximation of the Dirac operator. In this section, we consider $L^2(\nu)$ as the Hilbert space $\mathcal{H}$.
\begin{mdef}[Semigroup] A \textit{one-parameter unitary group} is a map $t\rightarrow P_t$ from $\mathbb{R}_+$ to $\mathcal{L}(\mathbb{H})$ such that 
\begin{equation}
P_0=1 \qquad P_{t+s}=P_tP_s,
\end{equation}
and $t\rightarrow P_t$ is continuous in the strong topology, i.e. $U_t\xrightarrow[]{s}U_{t_0}$ when $t\rightarrow t_0$.
\end{mdef}
Given a semigroup $P_t$ in $L^2$, define the \textit{generator} $\mathscr{L}$ of the semigroup by
\begin{equation}
\mathscr{L}(f):=\lim_{t\rightarrow 0} \frac{f-P_tf}{t},
\end{equation}
where the limit is understood in the $L^2$-norm. The \textit{domain} $\text{dom}(\mathscr{L})$ of the generator $\mathscr{L}$ is the space of functions $f\in \mathcal{H}$ for which the above limit exists. By the Hille-Yosida theorem, $\text{dom}(\mathscr{L})$ is dense in $L^2$. Moreover, $P_t$ can be recovered from $\mathscr{L}$ as follows:
\begin{equation}
P_t=\exp(-t\mathscr{L}).
\end{equation}
understood in the sense of spectral theory.\\
We then consider the operator $L=-i\frac{d}{dx}$ on $\mathcal{H}$ with $\text{dom}(\mathscr{L})=\lbrace f\in L^2(\mathbb{R}) : \xi\widehat{f}\in L^2(\mathbb{R})\rbrace $. Recall that $L$ is unitary equivalent to the left-multiplication operator $M_\xi$ using the Fourier transform
\begin{equation}
\mathcal{F}L\mathcal{F}^{-1}=\xi \widehat{f}.
\end{equation}
Then the associated semigroup $U_t$, so-called \textit{momentum operator}, is given by the left-multiplication operator in Fourier basis: $\mathcal{F}U_t\mathcal{F}^{-1}\widehat{f}=\xi \widehat{f}$. Therefore,
\begin{equation}
U_tf(x)=\mathcal{F}^{-1}(e^{it\xi}\widehat{f})(x)=\int_{\mathbb{R}}e^{i(x+t)\xi}\widehat{f}=f(x+t).
\end{equation}
We can then use $U_t$ in the definition of the Dirac operator $D$. This is what we have done to some extend (although not presented in the following framework) in our previous paper \cite{tageddine_noncommutative_2022}.
\section{Integral operators on manifolds} 
In this section, we extend the previous results to the case where $M$ is a smooth manifold of dimension $d$. Since we are interested in the Dirac operator over \textit{spin manifolds}, these results are crucial for the rest of the present work.
\subsection{Wrapped distributions}
Let us recall that on a well-suited open neighbourhood of a spin manifold $M$, the Dirac operator can be written as in Equation \eqref{DiracOp}. In the previous section, we have exhibited a probability distribution:
\begin{equation}
T_{\mu_{y,t}}(f)(y,t)=\int_{\mathbb{R}^d}\varphi(x)G_{t}(x,y)dx
\end{equation}
that we are going to use to approximate the partial derivatives appearing in the expression of the Dirac using the relation:
\begin{equation}
\partial_t|_{t=0}f=\partial_if
\label{FP}
\end{equation}
satisfied by the function $f$ obtained using the Green function $G$. In order to do so, we will need to extend the distribution $T_\mu$ defined on $\mathbb{R}^d$ to a compact manifold $M$. This will be done using \textit{wrapped distributions}.\par
Given the Lebesgue measure on $T_xM$ (in this section the tangent space $T_xM$ is identified with $\mathbb{R}^d$), the Riemannian volume form $\text{vol}_g$ on $M$ and a diffeomorphism $\phi:M \supset V\rightarrow U$, one can define a density $h$ on $M$ (whose support is included in $V\subset M$) to a density on $U\subset T_xM$. This is done by using $\phi$ as a push-forward map. The construction goes as follows: given a volume form $\omega$ written in coordinates as
\begin{equation}
\omega_g=h(x)\text{vol}_g
\end{equation}
Then, the integration on $V$ of this volume form is given by:
\begin{equation}
\int_V\omega_g:=\int_Uh(\phi(x)) |\det(d_x\phi)|dx_1\cdots dx_n
\end{equation}
where the integral is written using $\phi$ as a coordinate chart.\par
Among choices of $\phi$ an interesting candidate is the exponential map at point $p\in M$, $\exp_p:T_pM\rightarrow M$, due to its algebraic and geometric properties.
\begin{mprop} Let $(M,g)$ be a Riemannian manifold. Fore every point $p\in M$, there is an open subset $W\subseteq M$, with $p\in W$ and a number $\epsilon >0$, so that:
\begin{equation}
\exp_q:B(0,\epsilon)\subseteq T_qM\rightarrow U_q=\exp(B(0,\epsilon))\subseteq M
\end{equation}
is a diffeomorphism for every $q\in W$, with $W\subseteq U_q$.
\end{mprop}
\begin{mdef}[Normal neighbourhood] Let $(M,g)$ be a Riemannian manifold. For any $q\in M$, an open neighbourhood of $q$ of the form $U_q=\exp_q(B(0,\epsilon))$ where $\exp_q$ is a diffeomorphism from the open ball $B(0,\epsilon)$ onto $U_q$, is called a \textit{normal neighbourhood}.
\end{mdef}
\begin{mdef}[Injectivity radius] Let $(M,g)$ be a Riemannian manifold. For every point $p\in M$, the \textit{injectivity radius of $M$ at $p$}, denoted $\delta(p)$, is the least upper bound of the numbers $r>0$, such that $\exp_p$ is a diffeomorphism on the open ball $B(0,r)\subseteq T_pM$. The \textit{injectivity radius, $\delta(M)$ of $M$} is defined as: 
\begin{equation}
\delta(M):=\inf_{p\in M}\delta (p).
\end{equation}
\end{mdef}
In what will follow, we will simply denote by $\delta$ the \textit{injectivity radius} of $M$; we will also restrict to manifolds with strictly positive injectivity radius. Let $p\in M$, we then consider the exponential map $\exp_p:B(0,\delta)\rightarrow \exp_p(B(0,\delta))$. We associate to the density defined in Equation \eqref{density}, the volume form
\begin{equation}
\omega_{0,t}=\frac{1}{(4\pi \Phi(t))^{\frac{d}{2}}}\exp\left(-\frac{|\exp_p^{-1}(x)+st|^2}{4\Phi(t)}\right)\text{vol}_g=G_t(\exp_p^{-1}(x),0)\text{vol}_g
\end{equation}
We use the inverse of the exponential map to pushforward the Green function from the tangent space $T_pM$ at $p$ (identified with $\mathbb{R}^d$) to the manifold $M$. Integration on the domain of injectivity given by the injectivity radius, we obtain 
\begin{equation}
\int_{\exp_p(B_\delta)}f\omega_{0,t} = \int_{B_\delta}f\circ \exp_p(x)G_{0,t}(x)\det(d_{x}\exp_p)dx
\end{equation}
If we denote by $\left\lbrace e_1,\dots,e_d \right\rbrace $ a local orthogonal frame in $TM$, then by taking $s=e_i(p)$, we have proven the following lemma.
\begin{mlem} 
\label{Pushfor}
Consider the tangent vector $e_i(p)\in T_pM$ for $i\in \left\lbrace 1,\dots, d\right\rbrace $ and let $s=e_i(p)$. Then, the following holds:
\begin{equation}
\int_{\exp_p(B_\delta)}f\omega_{0,t}=\frac{1}{(4\pi \Phi(t))^{\frac{d}{2}}}\int_{B_\delta}\widetilde{f}(x)\exp\left(-\frac{|x+e_i(p)t|^2}{4\Phi(t)}\right)\det(d_{x}\exp_p)dx
\end{equation}
where $\exp_p$ is the exponential map on the Riemannian manifold $(M,g)$ and the function $\widetilde{f}(x)=f(\exp_p(x))$.
\end{mlem}
More generally, we can define the following map on the whole manifold $M$:
\begin{equation*}
\Theta : M\rightarrow C_0(\mathbb{R}_+),\quad p\mapsto \left\lbrace t\mapsto T_{\omega_{0,t}}(f) = \int_{\exp_p(B_\delta)}f\omega_{0,t}\right\rbrace_{t>0}.
\end{equation*}
This map associates to any normal neighbourhood a function of the variable $t$. We are going to show that each family of operators satisfies Equation \eqref{FP}.
\begin{mdef}[Jacobi field \cite{carmo_riemannian_1992}] Let $p\in M$ and $\gamma:\left[ 0,a\right] \rightarrow M$ be a geodesic with $\gamma(0)=p,\gamma'(0)=v.$ Let $w\in T_v(T_pM)$ with $|w|=1$. A Jacobi field $J$ along $\gamma$ given by 
\begin{equation}
J(t)=(d\exp_p)_{tv}(tw).
\end{equation}
\end{mdef}
\begin{mlem} Let $J$ be a Jacobi field. We have the following Taylor expansion about $t=0$:
\begin{equation}
\left\langle w, J(t)\right\rangle  = t + r(t),
\end{equation}
where $\lim_{t\rightarrow 0}\frac{r(t)}{t^2}=0$.
\label{Jacob}
\end{mlem}
\begin{proof}
From the definition of $J$ and the properties of the exponential map, we have that $J(0)=(d_0\exp_p)(0)=0$ and $J'(0)=w$. Hence, the first two coefficients of the Taylor expansion are
\begin{align*}
\left\langle w, J(0) \right\rangle =0,\\
\left\langle w, J'(0) \right\rangle = 1.
\end{align*}
As $J$ is a Jacobi field we have $J''(0)=-R(\gamma',J(0))\gamma'(0)=0$, where $R$ is the curvature tensor. This yields, 
\begin{equation}
\left\langle w, J''(0) \right\rangle = 0,
\end{equation}
which concludes the proof.
\end{proof}
\begin{mlem}
\label{taylor}
Define the smooth map:
\begin{equation}
G:T_pM\rightarrow \mathbb{R}, \qquad y\mapsto \det \left( d_y\exp_p\right),
\end{equation}
then, it satisfies $\nabla(G)(0)=0$.
\end{mlem}
\begin{proof}
In order to compute $\nabla (G)(0)$, we first use Jacobi's identity 
\begin{equation}
\left. \frac{d}{dt}\det\left( d_{ty}\exp_p\right)\right|_{t=0}=\det(d_{0}\exp_p)\text{tr}\left(d_{0}\exp_p^{-1}\left. \frac{d}{dt}\right|_{t=0}d_{ty}\exp_p\right) 
\end{equation}
which simplifies into
\begin{equation}
\left. \frac{d}{dt}\det\left( d_{ty}\exp_p\right)\right|_{t=0}=\text{tr}\left(\left. \frac{d}{dt}\right|_{t=0}d_{ty}\exp_p\right).
\end{equation}
Using the definition of a Jacobi field and linearity of $\text{tr}$, we have that
\begin{align*}
\text{tr}\left(\left. \frac{d}{dt}\right|_{t=0}d_{ty}\exp_p\right) &= \sum_{i=1}^d\left\langle v_i, \left( \left. \frac{d}{dt}\right|_{t=0}d_{ty}\exp_p\right) v_i\right\rangle, \\
&=\left. \frac{d}{dt}\right|_{t=0}\sum_{i=1}^d\left\langle v_i, d_{ty}\exp_p(v_i)\right\rangle,\\
&=\left. \frac{d}{dt}\right|_{t=0}\sum_{i=1}^d\frac{1}{t}\left\langle v_i,J_i(t)\right\rangle.
\end{align*}
Now using the Taylor expansion obtained in Lemma \eqref{Jacob}, we get:
\begin{equation}
\left\langle v_i,J_i(t)\right\rangle = t + r(t)
\end{equation}
where $r(t)=o(t^2)$, we conclude that:
\begin{equation}
\text{tr}\left(\left. \frac{d}{dt}\right|_{t=0}d_{ty}\exp_p\right) = 0.
\end{equation}
\end{proof}
\begin{mth}
\label{limit}
The following limit holds at $p\in M$
\begin{equation}
\left. \frac{\partial}{\partial t}\left( \int_{\exp_p(B_\delta)}f\omega_{0,t}\right)\right|_{t=0}=e_i(f)(p).
\end{equation}
\label{Result}
\end{mth}
\begin{proof} 
The result follows from the property of the Green function given in Lemma \ref{FP-th} and the reduction to an open ball obtained in Lemma \label{Reduc}. Then, using the integration equality given in Lemma \ref{Pushfor} on the open ball $B_\delta$ and the isomorphism $T_0(T_p(M))\simeq T_p(M)$, we have:
\begin{equation*}
\left. \frac{\partial}{\partial t}\left( \int_{\exp_p(B_\delta)}f\omega_{0,t}\right)\right|_{t=0}=e_i\left( \exp_{p*}(f)\det \left( d\exp_p\right)\right)(0).
\end{equation*}
Finally, to conclude we use Lemma \ref{taylor} and deduce that:
\begin{equation}
\left. \frac{\partial}{\partial t}\left(\int_{\exp_p(B_\delta)}f\omega_{0,t}\right)\right|_{t=0}=e_i(f)(p).
\end{equation}
\end{proof}
\section{Statistical fluctuations of differential structures}
\label{sect3}
We are now ready to state and prove Theorem \ref{main}. We keep the same notations as the previous sections: $M$ is a compact Riemannian manifold of dimension $d$; we consider a point $p\in M$ and a normal neighbourhood $U_p$ associated to it; we denote by $\left\lbrace e_1,\dots,e_d \right\rbrace $ a local orthogonal frame in $TM$. We finally define the orthonormal family of vectors $\left\lbrace s_1,\dots,s_d \right\rbrace $ such that:
\begin{equation}
s_j=e_j(p), \qquad \forall j\in \left\lbrace 1,\dots,d \right\rbrace .
\label{frame}
\end{equation}
\subsection{The Dirac operator}
We start by recalling that the notation $D_X$ means: a Dirac operator $D$ associated to a graph $X$ in the sense of \cite[Def. 4.3]{tageddine_noncommutative_2022}. Now, let $n$ be a positive integer and fix a graph $X_n$ equipped with a Dirac operator $D_{X_n}$ and with set of vertices $\lbrace x_1,\dots x_n\rbrace$. In addition, we are going to consider $n$ copies of the same graph $X_n$, each of which is equipped with a Dirac operator $D_{X_k}$ and with a set of vertices denoted by $\lbrace x^k_1,\dots x^k_n\rbrace$, for $1\leq k\leq n$. Then, we have a sequence of Dirac operators 
\begin{equation}
(D_{X_1},D_{X_2},\cdots, D_{X_n} ) \in \mathfrak{gl}_{2m_n}^-(\mathbb{C})^n,
\end{equation}
acting on a sequence of diagonal elements $(a_1, \cdots, a_n)$ with each $a_i\in\mathfrak{A}_n$. \\
If we denote by $(a_k^i)_{1\leq i \leq n}$ the coefficients of $a_k$ in the block diagonal, then using the projection maps $M\rightarrow X_k$ we can identify these values with evaluations of a smooth functions, denoted by $a$ (see \cite[Prop. 3.5]{tageddine_noncommutative_2022} for more details):
\begin{equation}
a_k^i=a(x^k_i), \quad \forall i \in \left\lbrace 1,\dots, n\right\rbrace 
\end{equation}
for some point $x^k_i\in M$. Fix a point $p\in M$ and a neighbourhood $U_p$ of $p$ in $M$. Then, consider a sequence of points $\left\lbrace x^k_1, \dots, x^k_n \right\rbrace $ in $U_p$, for $1\leq k\leq n$, such that, for a chosen index $i_0$ (not depending on $k$), we have $x^k_{i_0}=p$.
We then define the coefficients $(\omega^k_{ij})_{1\leq i,j\leq n}$ of $D_{X_n}$ as follows:
\begin{align*}
&\omega^k_{ij}(\hbar)=\frac{1}{(4\pi \Phi(\hbar))^{\frac{d}{2}}}\exp\left(-\frac{|y^k_i+s_j\hbar|^2}{4\Phi(\hbar)}\right) \quad \text{for $1\leq i,j\leq n$ and for $1\leq k\le n$},\\
&\text{with $y_i^k:=\exp^{-1}_p(x_i^k)$}.
\end{align*}
Furthermore, for every integer $1\leq k \leq n$, we define a family of projection elements such that $e_k\in M_{2m_n}(\mathbb{C})$ and we have the following matrix form:
\begin{equation}
e_kD_{X_k}e_k^*=\begin{pmatrix}
\bigzero 
  & \rvline & 
   \begin{array}{ccccc}
 &  &  &  &   \\ 
 &  & \bigzero  &  &   \\ 
* & *  & \omega^k_{i_{0}j} & *  & *  \\ 
 &  &  &  &   \\ 
 &  & \bigzero &  &   
\end{array} \\
\hline
  \begin{array}{ccccc}
 &  & * &  &   \\ 
 &  & * &  &   \\ 
 & \bigzero   & \omega^k_{i_{0}j} &  \bigzero  &   \\ 
 &  & * &  &   \\ 
 &  & * &  &   
\end{array} 
  & \rvline & \bigzero 
\end{pmatrix}.
\label{Diracbis}
\end{equation}
\begin{mrmk}The non-zero coefficients correspond to the adjacency points of $i_0$.
\end{mrmk}
Hence, if we recall the expression given by the commutator in Equation \eqref{dirac}, we consider the following average of operators over the $n$ copies of $X_n$:
\begin{equation}
\widehat{S}_n^{\hbar_n} (a) :=\frac{1}{n}\sum_{k=1}^ne_k\left[ D_{X_k}, a_k\right]e_k^* =\frac{i}{n\hbar_n}\sum_{k=1}^n\sum_{j=1}^{d+1}\omega^k_{i_0j}(\hbar_n)\alpha_{i_0j}(a_k)Y\otimes E_{i_0j},
\end{equation}
where $\alpha_{i_0j}(a_k)=a(x^k_j)-a(p)$. Moreover, for the purpose of the proof of the main theorem, we define a second operator given by: 
\begin{equation}
S_{j,n}^{\hbar_n}:C^\infty(M)\rightarrow \mathbb{R}, \qquad S_{j,n}^{\hbar_n}(a)=\frac{1}{n\hbar}\sum_{k=1}^n\omega^k_{i_0j}(\hbar_n)\alpha_{i_0j}(a_k).
\label{averop}
\end{equation}
We assume now that the points $\left\lbrace x^k_1, \dots, x^k_n \right\rbrace $ are thought as random variables independent and identically distributed (i.i.d.) from a uniform distribution. Let us recall the definition of the map $\Psi$ given in Equation \eqref{mapPsi}:
\begin{equation}
\Psi: V_{i_0}\rightarrow Cl(\mathbb{R}^d), \quad \Psi\left( \left[ D,a\right]_{i_0}\right) =  \frac{i}{\hbar}\sum_{j=1}^{d}\omega_{i_0j}\alpha_{i_0j}(a)e_j.
\end{equation}
Then, we can prove the following theorem.
\begin{mth}Let $\left\lbrace x^k_{i_0}\right\rbrace _{k=1}^n$ be a sequence of i.i.d. sampled points from a uniform distribution on a open normal neighbourhood $U_p$ of a point $p$ in a compact Riemannian manifold $M$ of dimension $d$. Let $\tilde{S}^{\hbar_n}_n$ be the associated operator given by:
\begin{equation}
\widehat{S}_n^{\hbar_n}:C^\infty(U_p)\rightarrow M_2(\mathbb{R})\otimes U(\mathfrak{gl}_{2m_n}), \qquad \widehat{S}_n^{\hbar_n} (a) :=\frac{1}{n}\sum_{k=1}^ne_k\left[ D^k_X, a_k\right]e_k^*.
\end{equation}
Put $\hbar_n=n^{-\alpha}$, where $\alpha>0$, then for $a\in C^\infty(U_p)$, in probability:
\begin{equation*}
\lim_{n\rightarrow\infty}\Psi\circ\widehat{S}_n^{\hbar_n}(a)=\left[ \mathcal{D},a\right](p).
\end{equation*} 
\end{mth}
\begin{proof}
We consider the average operator defined by Equation \eqref{averop}. It is then sufficient to prove that for $\hbar_n=n^{-\alpha}$, where $\alpha>0$, and for $a\in C^\infty(U_p)$, we have:
\begin{equation*}
\lim_{n\rightarrow\infty}S_{j,n}^{\hbar_n}(a)=e_j(a)(p) \quad \forall 1\leq j\leq d.
\end{equation*} 
in probability, and then apply the map $\widehat{\psi}$ defined in \eqref{mapx}. Recall that $e_j$ is given in Equation \eqref{frame}.\\

The expectation value of the random variable $S_n^{\hbar_n}(a)$ is given by:
\begin{equation}
\mathbb{E}S_{j,n}^{\hbar_n} (a)(p)=\frac{1}{\hbar_n}\int_{\exp_p(B_\delta)}\omega_{0,\hbar_n}(a-a(p)),
\end{equation}
where we assume, without lost of generality, that the volume of $M$ is equal to one. We recognize then an approximation of the time derivative at $0$ in Equation \eqref{Result}. Thus, applying Hoeffding's inequality, we have:
\begin{equation}
\mathbb{P} \Big[ \left|S^{\hbar_n}_{j,n}(a)(p)-\mathbb{E}S_{j,n}^{\hbar_n}(a)(p)\right|>\varepsilon \Big] \leq 2\exp\left(-\frac{\varepsilon^2n}{KC_d(n^\alpha)^2}\right).
\end{equation}
Choosing $\hbar$ as a function of $n$, such that $\hbar(n)=n^{-\alpha}$, where $\alpha>0$, we have, for any real number $\varepsilon>0$:
\begin{equation}
\lim_{n\rightarrow \infty}\mathbb{P} \Big[  \left|S^{\hbar_n}_{j,n}(a)(p)-\mathbb{E}S_{j,n}^{\hbar_n}(a)(p)\right|>\varepsilon \Big] = 0.
\end{equation}
Finally, we prove the statement using Theorem \ref{limit}:
\begin{equation}
\lim_{n\rightarrow \infty}S_{j,n}^{\hbar_n}(a)(p)=e_j(a)(p),
\end{equation}
along with the definition of the map $\Psi$ in Equation \eqref{mapPsi}.
\end{proof}
\begin{mrmk}
Every line in the matrix $D_X$ corresponds then to a point $p$ and a normal neighbourhood $U_{p}$ obtained from the image of the exponential map of a ball of radius $\delta$. Indeed, since $M$ is compact, we have a finite cover $\lbrace U_{p_i}\rbrace_{i=1}^N$ with centre $\lbrace p_i\rbrace_{i=1}^N$ every one of which being associated to a line of $D_X$. \\
If we consider a sequence of Dirac operators $D_{X_n}$, then what we are doing is in fact taking refinements of normal neighbourhoods, increasing with the numbers of vertices in $X_n$.
\end{mrmk}
\subsection{Uniform convergence}
It is interesting to mention that the previous result can be extended to have a uniform convergence, following the same steps as \cite[Prop. 6.1]{belkin_towards_nodate}. We then state the result without proof.
\begin{mprop}
Let $\mathscr{F}$ be an equicontinuous family of functions with a uniform bound up to the second derivative. Then for each $\hbar >0$, we have:
\begin{equation}
\lim_{n\rightarrow\infty}\mathbb{P}\left[\sup_{a\in\mathscr{F}} \ \Big |S^{\hbar_n}_n(a)(p)-\mathbb{E}S_n^{\hbar_n}(a)(p)\Big |>\varepsilon \right] = 0. 
\end{equation}
\end{mprop}
\begin{mth}
\label{mainDirac}
Let $\mathscr{F}$ be a family of smooth functions with uniformly bounded derivatives up to the second order. Let $\left\lbrace x^k_{i_0}\right\rbrace _{k=1}^n$ be a sequence of i.i.d. sampled points from a uniform distribution on an open normal neighbourhood $U_p$ of a point $p$ in a compact Riemannian manifold $M$ of dimension $d$. Let $\widehat{S}_n^{\hbar_n}$ be the associated operator given by:
\begin{equation}
\widehat{S}_n^{\hbar_n}:C^\infty(U_p)\rightarrow M_2(\mathbb{R})\otimes U(\mathfrak{gl}_{2m_n}), \qquad \widehat{S}_n^{\hbar_n} (a) :=\frac{1}{n}\sum_{k=1}^ne_k\left[ D^k_X, a_k\right]e_k^*.
\end{equation}
 Put $\hbar_n=n^{-\alpha}$, where $\alpha>0$, then in probability:
\begin{equation}
\lim_{n\rightarrow\infty} \sup_{a\in\mathscr{F}} \ \Big |\Psi\circ\widehat{S}^{\hbar_n}_n(a)(p)-\left[ \mathcal{D},a\right](p)\Big |=0.
\end{equation}
\end{mth}
\subsection{The Laplacian }
In this final section, we want to study the convergence result for the Laplacian defined by Equation \eqref{Laplacianop}. If we take the second derivative in time in the initial value problem \eqref{HJeqt}, we see that
\begin{equation}
\partial_t^2u=s\cdot \nabla \partial_t u +\varphi\Delta\partial_tu+\varphi'\Delta u 
\end{equation}
from which we see that if $\varphi(t) =t$ then when we take the limit when $t$ goes to zero and obtain:
\begin{equation}
\partial_t^2u|_{t=0}=(s\cdot \nabla)^2u_0+\Delta u_0.
\end{equation}
Moreover, if we let $\widetilde{G}(x,t)$ defined by
\begin{equation}
\widetilde{G}(x,t) = \sum_{j=1}^m \lambda_jG_j(x,t), \qquad \sum_{\lambda_j=1}^m\lambda_j=1
\end{equation}
where $G_j(x,t)$ is the Green function of Equation \eqref{HJeqt} with unit vector $s_j$ and the $\lambda_j$ are here to ensure that $\widetilde{G}$ remains a probability distribution. Then, by linearity the function 
\begin{equation}
\widetilde{u}=u_0*\widetilde{G}
\end{equation}
satisfies the equation:
\begin{equation}
\partial_t^2\widetilde{u}|_{t=0}=(\tilde{s}\cdot \nabla)^2u_0+\Delta u_0, \qquad \text{with $\widetilde{s}=\sum_{j}s_j$}.
\end{equation}
Hence, if we pick the vectors $s_j$ such that $\widetilde{s}$ is zero, then we are left with the following equations:
\begin{equation}
\partial_t^2\widetilde{u}|_{t=0}=\Delta u_0, \quad \text{and} \quad \partial_t\widetilde{u}|_{t=0}=0,
\end{equation}
where, the second equation is obtained from Equation \eqref{HJeqt} after taking $t$ to zero. Therefore, in $\mathbb{R}^d$, we have for any smooth initial condition $f$, the following limit:
\begin{equation}
    \lim_{t\rightarrow 0}\frac{1}{t^2}\int_{\mathbb{R}^d}\widetilde{G}(x,t)(f(x)-f(0))dx=\Delta f(0).
\end{equation}
In addition, if the dimension $d$ is greater than two, then the reduction  to a ball Lemma \ref{Reduc} still holds. In the rest of this section, we will assume that $d\geq 2$. Hence, we can extend the distribution $\widetilde{G}$ to a manifold as in Section \ref{sect3}. 
The volume form $\omega_{0,t}$ is now given by:
\begin{equation}
\omega_{0,t}=\frac{1}{(4\pi \Phi(t))^{\frac{d}{2}}}\sum_{j=1}^{d+1}\lambda_j\exp\left(-\frac{|\exp_p^{-1}(x)+s_jt|^2}{4\Phi(t)}\right)\text{vol}_g.
\end{equation}
We recognize a convex combination of the Green function obtained in the previous section. Consequently, following the same steps as in Theorem \ref{FP-th}, we see that this distribution satisfies the equation:
\begin{equation}
\left. \frac{\partial^2}{\partial t^2}\left( \int_{\exp_p(B_\delta)}\omega_{0,t}f\right)\right|_{t=0}=\Delta_M (f)(p).
\end{equation}
 Now, the Laplacian $\Delta_{X_k}$ obtained from the Dirac operator \eqref{Diracbis} and acting on an element $a_k$ is given by:
\begin{equation}
\Delta_{X_k}(a_k)=\frac{1}{\hbar^2}\sum_{j=1}^{d+1}(\omega^k_{ij})^2\alpha_{ij}(a_k)J.
\end{equation}
Then, we assume that the coefficients $\omega_{ij}^k$ of the Dirac $D_{X_k}$ are given by:
\begin{equation}
\omega^k_{ij}(\hbar)=\frac{1}{(4\pi \hbar^2)^{\frac{d}{4}}}\sqrt{\lambda_j}\exp\left(-\frac{|y^k_i+s_j\hbar|^2}{8\hbar^2}\right),
\end{equation}
where $\lambda_j$ are positive numbers to be specified. Therefore, we are lead to study the convergence of the averaging operator:
\begin{equation}
\Omega_n^{\hbar}(a)(p) = \frac{1}{(4\pi \hbar^2)^{\frac{d}{2}}n\hbar^2}\sum_{k=1}^n\sum_{j=1}^{d+1}\lambda_j\exp\left(-\frac{|y^k_i+s_j\hbar|^2}{4\hbar^2}\right)\alpha_{ij}(a_k).
\end{equation}
such that, if $u=\sum_{i=1}^{d}s_i$, then $s_{d+1}=-u/\|u\|$ and $\lambda_j=1/(d+\|u\|)$ for $1\leq j\leq d$ and $\lambda_{d+1}=\|u\|/(d+\|u\|)$. Notice then that $\sum_{j=1}^{d+1}\lambda_j=1$.\\
Moreover, the expectation value of the random variable $\Omega_n^{\hbar_n}(a)$ is given by:
\begin{equation}
\mathbb{E}\Omega_n^{\hbar}(a)(p) =\int_{\exp_p(B_\delta)}\omega_{0,\hbar_n}\left( a(x)-a(p)\right)
\end{equation}
\begin{mth}
\label{laplace}
Let $\mathscr{F}$ be a family of smooth functions with uniformly bounded derivatives up to the third order. Let $\left\lbrace x_i\right\rbrace _{i=1}^n$ be a sequence of i.i.d. sampled points from a uniform distribution on an open normal neighbourhood $U_p$ of a point $p$ in a compact Riemannian manifold $M$ of dimension $d$. $\Omega_n^{\hbar_n}:C^\infty(U_p)\rightarrow \mathbb{R}$ be the associated operator given by:
\begin{equation*}
 \Omega_n^{\hbar_n}(a)(p) = \frac{1}{(4\pi \Phi(\hbar))^{\frac{d}{2}}n\hbar^2}\sum_{k=1}^n\sum_{j=1}^{d+1}\lambda_jexp\left(-\frac{|y_i^k+s_jt|^2}{4\Phi(t)}\right)\alpha_{ij}(a_k).
\end{equation*}
 Put $\hbar_n=n^{-\alpha}$, where $\alpha>0$, then in probability:
\begin{equation}
\lim_{n\rightarrow\infty} \sup_{a\in\mathscr{F}} \ \Big |\Omega_n^{\hbar_n}(a)(p)-\Delta_M (a)(p)\Big |=0
\end{equation}
\end{mth}
\subsection{Discussion} 
Going back to the definition of the Dirac operator associated to a graph $X$ with non-zero coefficients $\omega_{ij}$, we recall that the goal was to compute the values $\omega_{ij}$ in order to obtain a convergence when considering a sequence of refined triangulations. We have exhibited the coefficients 
\begin{equation}
\omega_{ij}(\hbar)=\frac{1}{(4\pi \Phi(\hbar))^{\frac{d}{2}}}\exp\left(-\frac{|y_i+s_j\hbar|^2}{4\Phi(\hbar)}\right) 
\end{equation}
obtained from the Green function given in Equation \eqref{Green}. Hence, we are able to prove a convergence result to the Dirac operator on a normal neighbourhood (Theorem \ref{mainDirac}) as well as a convergence of the Laplace operator (Theorem \ref{laplace}). However, as far as the Laplacian is concerned, this choice is not unique, in fact one could take the values of $\omega_{ij}$ obtained from a normal distribution and such that:
\begin{equation}
\omega^2_{ij}(\hbar)= \exp\left(\frac{\|x_i-x_j\|^2}{4\hbar} \right)
\end{equation}
and still get a convergence result. Nevertheless, keeping in mind that we are also interested in the convergence of the square root i.e. to the Dirac operator, it is not clear that such a choice of coefficients would also work.\par

Moreover, one may also consider classical discretizations of the Laplacian such as the combinatorial one with the choice:
\begin{equation}
\text{vertices $i$ and $j$ do not share an edge} \Leftrightarrow \omega_{ij}=0, \ \forall i,j
\end{equation}
or the \textit{cotangent Laplacian} with the choice
\begin{equation}
\omega_{ij}^2=\left\lbrace 
\begin{array}{cc}
\frac{1}{2}\left(\cot \alpha_{ij} +\cot \beta_{ij} \right)  & \text{$ij$ is an edge}, \\ 
-\sum_{k\sim i}\omega^2_{ik} & i=j, \\ 
0 & \text{otherwise}.
\end{array} 
\right. 
\end{equation}
In these two cases the question of convergence to the Laplacian is unclear \cite{bobenko_convergence_2008}, let alone convergence in the square root.\par
There is therefore an important direction worth investigating: whether a convergent Laplacian constructed from a specific distribution or obtained from a known discretization implies convergence of its associated Dirac operator.

\bibliographystyle{amsplain}
\bibliography{Article2Bibli.bib}
\end{document}